# CHARACTERIZING AUTOMORPHISM GROUPS OF ORDERED ABELIAN GROUPS

Rüdiger Göbel and Saharon Shelah

In this short note we want to characterize the groups isomorphic to full automorphism groups of ordered abelian groups. The result will follow from classical theorems on ordered groups adding an argument from proofs used to realize rings as endomorphism rings of abelian groups, see [2]. Recall that $H$ is a right ordered group (RO-group) if $(H, \cdot)$ is a group and $(H, <)$ is a linear order satisfying the following compatibility condition

(RO)  For all $h < g, k$ in $H$ follows $hk < gk$.

Similary define (LO), the left compatibility condition. If $(H, \cdot, <)$ satisfies both (RO) and (LO) then $H$ is an ordered group. Obviously abelian RO-groups are ordered groups, in which case we often replace multiplication by $+$. A group $H$ is right orderable if it permits a linear order which makes it an RO-group. We do not distinguish RO-groups and groups which are right orderable. From the fact that cyclic ordered groups are infinite, it is clear that RO-groups are torsion-free. By an old theorem of Smirnov a group is an RO-group if and only if it is (isomorphic to) a subgroup of $\operatorname{Aut}(A, +, <)$ of an ordered free abelian group $A$, see Mura, Rhemtulla [5, p. 129, Theorem 7.1.3]. We will use the obvious representation as subgroup of $\operatorname{Aut}(A, +, <)$ below. On the other hand there are torsion-free groups, in fact polycyclic groups, which are not RO-groups, a result due to Smirnov, see [5, p. 127]. Note that torsion-free polycyclic groups are even finitely generated - iterated extensions of $\mathbb{Z}$. Our main result then reads as follows.

**Theorem 1** *For a group $H$ the following are equivalent.*

*(1) $H$ is an RO-group.*

*(2) There is an ordered abelian group $G = (G, +, <)$ with $\operatorname{Aut}(G, +, <) \cong H$.*

---

[0][GbSh:780] in Shelah's list of publications

(3) If $K$ is any ordered field, then we find an ordered extension field $F$ such that $\mathrm{Aut}\,(F) \cong H$.

This result is in sharp contrast to Corner's [1] result classifying all finite groups which are automorphisms groups of abelian groups, where many groups like $\mathbb{Z}/7\mathbb{Z}$ do not come up. The equivalence $(1), (3)$ is taken from Dugas, Göbel [3], and (1) follows from (3) by the above theorem of Smirnov. It remains to show that (1) implies (2), in fact we will show a stronger implication

($2^*$) There is an $\aleph_1$-free ordered abelian group $G = (G, +, <)$ with $\mathrm{Aut}\,(G, +, <) \cong H$.

Here $G$ is $\aleph_1$-free if all its countable subgroups are free. Consider the group ring $R = \mathbb{Z}H$ and let $B = \bigoplus R$ be a "large enough" free $R$-module. Hence $R \subseteq \mathrm{End}\,_R B$ by scalar multiplication with elements from $R$ on the right of $B$. We will construct an $R$-module $G$ such that $B \subseteq_* G \subseteq_* \widehat{B}$ and $\mathrm{End}\,A = R$. Here $\widehat{B}$ is the $S$-adic completion of $B$ with respect to some suitable, multiplicatively closed subset $S \subseteq \mathbb{N} \subseteq R$; e.g. $S = \{p^n \mid n \in \omega\}$. Moreover "$\subseteq_*$" denotes an $S$-pure submodule. It will be important that

$$G = \bigcup_{\alpha \in \lambda^*} G_*$$

is the union of an ascending continuous chain of $S$-pure $R$-submodules

$$G_\alpha \subseteq_* \widehat{B} \text{ with } G_0 = B \text{ and } G_{\alpha+1} = \langle G_\alpha, g_\alpha R \rangle_*$$

such that $\mathrm{Ann}\,_R g_\alpha = 0, \ G_\alpha \cap g_\alpha R = 0$

and either $G_{\alpha+1}/G_\alpha \cong R$ or $G_{\alpha+1}/G_\alpha \cong S^{-1}R$ is $S$-divisible.

Here $\langle G_\alpha, g_\alpha R \rangle_*$ denotes the smallest subgroup of $\widehat{B}$ which is $S$-pure and contains $G_\alpha, g_\alpha R$. This part will follow by arguments we used in several earlier papers, e.g. in [2]. As $H$ is assumed to be an RO-group by (1) we will turn $R = \mathbb{Z}H$ into a linear order satisfying the compatibility condition (RO) for multiplication with positive elements in the group ring.

**Proposition 2** *If $H$ is an RO-group, then the group ring $R = \mathbb{Z}H$ has a natural linear order satisfying (RO) for multiplication with positive elements. The monoid of positive elements of $R$ will be denoted by $R^{>0}$.*

We postpone the proof of Proposition 2 and assume it holds. It follows that $(R^+, <)$ is an ordered free abelian group and $R^{>0} \subseteq \mathrm{End}\,(R, +, <)$ thus also

$$R^{>0} \subseteq \mathrm{End}\,(B, +, <)$$

where the linear order will be extended and $B$ becomes an ordered abelian group.

Inductively we want to extend the order on to $G$. If $P_\alpha = \{g \in G_\alpha, 0 < g\}$ denotes the positive cone of $G_\alpha$, then we want to define the positive cone $P_{\alpha+1}$ of $G_{\alpha+1} = \langle G_\alpha, g_\alpha R \rangle_*$.

If $y \in G_{\alpha+1}$, there is $s \in S$ such that $ys = x + g_\alpha r$ for some $r \in R$ and $x \in G_\alpha$. Thus we define

$$y \in P_{\alpha+1} \iff r > 0 \text{ or } r = 0 \text{ and } x \in P_\alpha.$$

It is easy to see that $G_{\alpha+1} = -P_{\alpha+1} \cup P_{\alpha+1} \cup \{0\}$, and $P_{\alpha+1}$ is well-defined. If also $ys' = x' + g_\alpha r'$ then $xs' + g_\alpha rs' = x's + g_\alpha r's$. Hence $xs' - x's = g_\alpha(r's - rs') \in G_\alpha \cap g_\alpha R = 0$. From $\mathrm{Ann}\,_R g_\alpha = 0$ follows $r's = rs'$ and $s, s' > 0$ implies $r' > 0$ iff $r > 0$. Moreover $r = 0$ iff $r' = 0$ and in this case $xs' = x's$ thus $x \in P_\alpha$ iff $x' \in P_\alpha$. Note that $g_\alpha > 0$ follows from $1 > 0$. If $r' \in R^{>0}$ and $y \in P_{\alpha+1}$ as above, then $yr's = xr' + g_\alpha rr'$, hence either $r = 0$ and $xr' \in P_\alpha$ by induction hypothesis or $rr' > 0$ by Proposition 2, so $yr' \in P_{\alpha+1}$ and

$$R^{>0} \subseteq \mathrm{End}\,(G_{\alpha+1}, +, <)$$

follows from $R^{>0} \subseteq \mathrm{End}\,(G_\alpha, +, <)$. At limit ordinals $\beta < \lambda^*$ we take unions, thus $P_\beta = \bigcup_{\alpha < \beta} P_\alpha$ and it follows that $(G, +, <)$ is an ordered abelian group with

$$R^{>0} \subseteq \mathrm{End}\,(G, +, <).$$

If $r < 0$ then the action of $r$ on $B$, hence multiplication on a summand $Re$ of $B$ shows that $0e < 1e$ turns into $re < 0e$. Together with $R = \mathrm{End}\,G$ we get

$$R^{>0} = \mathrm{End}\,(G, +, <)$$

We derive the following

**Theorem 3** *If $H$ is an RO-group, $\lambda$ is any cardinal with $|H| \leq \lambda$ and $R^{>0}$ is the monoid of positive elements of the group ring $R = \mathbb{Z}H$, then there is an $\aleph_1$-free, ordered abelian group $(G, +, <)$ of cardinality $\lambda^{\aleph_0}$ with $R^{>0} = \mathrm{End}\,(G, +, <)$.*

provided Proposition 2 is shown and the used presentation of $G = \bigcup_{\alpha < \lambda} G_\alpha$ follows.

First we want to establish the last claim and note that $R^+$ is a free abelian group, in particular $R^+$ is cotorsion-free ($\mathrm{Hom}\,(\widehat{\mathbb{Z}}, R) = 0$), which is needed to apply Theorem 6.3 in Corner, Göbel [2, p. 465]. We need a very special case of Theorem 6.3, putting $\mathfrak{N} = \{0\}, J_1 = J = \emptyset$. Thus $\mathrm{End}\,G = R$ is immediate. It is easy to check that $G$ is $\aleph_1$-free. The group $G$ is obtained by transfinite induction as

$$\bigcup_{\alpha < \lambda^*} G_\alpha = G \text{ over } \alpha < \lambda^* \text{ with } |\lambda^*| = \lambda^{\aleph_0}$$

by using a weak version of Shelah's Black Box (see Appendix of [2]). Conditions $(II_0)$ $(II_\mu)$ and $(III_{\alpha+1})$ show that $G_{\alpha+1}$ is of the right form (replacing A by R), and Lemma 3.4 in [2, p. 456] for $N_\alpha^k = 0$ together with $(III_\alpha)$ implies that $G_\alpha \cap g_\alpha R = 0$, and $\mathrm{Ann}\,_R g_\alpha = 0$.

It remains to show Proposition 2.

If $r \in R$, write $r = \sum_{h \in H} r_h\, h$ with $r_h \in \mathbb{Z}$, similarly $r' = \sum_{h \in H} r'_h h$. We say that

$$r < r' \iff \exists\, h^* \in H,\ r_{h^*} < r'_{h^*} \text{ and } r_h = r'_h \forall\ h > h^*.$$

Let $[r] = \{h : r_h \neq 0\}$, then the positive cone of $R$ is

$$R^{>0} = \{r \in R : \exists \text{ maximal } h^* \in [r] \text{ and } r_{h^*} > 0\}.$$

It is easy to check that this is a linear ordering on $R$. From $R^{>0} \cdot R^{>0} \subseteq R^{>0}$ follows that multiplication with elements from $R^{>0}$ satisfies $(RO)$, thus $R^{>0} \subseteq \mathrm{End}\,(R, +, <)$. The ordering extends naturally to direct sums, see e.g. Theorem 2.1.1 in [5] thus $R^{>0} \subseteq \mathrm{End}\,(B, +, <)$.

Like in case of polynomial rings $\mathbb{R}[x]$ we can show the following

**Proposition 4** *If $H$ is an RO-group, $R = \mathbb{Z}H$ is the group ring and $U(R)$ are its units, then $U(R) = \pm H$.*

**Proof.** If $r = \sum_{h \in H} r_h h, r' = \sum_{h \in H} r'_h h \in R$ are as above with $rr' = r'r = 1$, then the product of the maximal coefficients $r_{h^*}$ and $r'_{h'^*}$ must be 1. This is only possible if $h^* h'^* = 1$ and if all other coefficients are 0. It follows that $r = r_{h^*} h^*$ and $r' = r_{h^*}^{-1} h_*^{-1}$, and $r_{h^*}, r_{h'^*}$ are units of the coefficient ring $\mathbb{Z}$. Hence $r, r' \in \pm H$.

**Remark** also follows from a more general result by Strojnowski on unique product groups, see Karpilovsky [4, p. 272, Corollary 8.4.8].

From Proposition 4 follows that the units of $R^{>0}$ are $U(R^{>0}) = H$. From Theorem 3 and $\mathrm{Aut}\,(G,<) = U(\mathrm{End}\,(G,+,<))$ follows our main result which immediately implies "(1) $\longrightarrow$ (2)" of Theorem 1.

**Corollary 5** *If $H$ is an RO-group of cardinality $|H| \leq \lambda$, then there is an $\aleph_1$-free, ordered abelian group $(G,<)$ of cardinality $\lambda^{\aleph_0}$ with $\mathrm{Aut}\,(G,+,<) = H$.*

Rüdiger Göbel
Fachbereich 6, Mathematik und Informatik
Universität Essen, 45117 Essen, Germany
e–mail: R.Goebel@Uni-Essen.De
and
Saharon Shelah
Department of Mathematics
Hebrew University, Jerusalem, Israel
and Rutgers University, Newbrunswick, NJ, U.S.A
e-mail: Shelah@math.huji.ae.il